\documentclass[11pt,a4paper]{article}

\usepackage[english]{babel} 
	\usepackage[T1]{fontenc} 
\usepackage{lmodern} 

\usepackage[margin=25mm]{geometry}

	\usepackage{booktabs}

	\usepackage{natbib}

	\usepackage[fleqn]{amsmath}
	\usepackage{amssymb,amstext,amsfonts,mathrsfs}
	\usepackage{array}
	\usepackage{cancel}
	\usepackage{lipsum}
	\usepackage{mathtools}
	\usepackage{url}
	\usepackage{color}
	\usepackage{ifthen}
	\usepackage{amsmath,blkarray}
	\usepackage{diagbox}
	\usepackage{caption}
	\usepackage{siunitx,booktabs,threeparttable,caption}
	\usepackage{mathabx}
	\usepackage{amsmath}
    \usepackage{algorithm}
    \usepackage{algpseudocode}
    \usepackage{lineno}

    \PassOptionsToPackage{noend}{algpseudocode}
\usepackage{algpseudocode}

\errorcontextlines\maxdimen

\makeatletter
\newcommand*{\algrule}[1][\algorithmicindent]{\makebox[#1][l]{\hspace*{.5em}\vrule height .75\baselineskip depth .25\baselineskip}}%

\newcount\ALG@printindent@tempcnta
\def\ALG@printindent{%
    \ifnum \theALG@nested>0
        \ifx\ALG@text\ALG@x@notext
            \addvspace{-3pt}
        \else
            \unskip
            \ALG@printindent@tempcnta=1
            \loop
                \algrule[\csname ALG@ind@\the\ALG@printindent@tempcnta\endcsname]%
                \advance \ALG@printindent@tempcnta 1
            \ifnum \ALG@printindent@tempcnta<\numexpr\theALG@nested+1\relax
            \repeat
        \fi
    \fi
    }%
\usepackage{etoolbox}
\patchcmd{\ALG@doentity}{\noindent\hskip\ALG@tlm}{\ALG@printindent}{}{\errmessage{failed to patch}}
\makeatother

\usepackage{tocloft} 

\newcommand{\problem}{RSPMST}

\allowdisplaybreaks

\usepackage{authblk}

\date{}

\title{Robust Optimization Approaches for Routing and Scheduling of Multi-Skilled Teams under Uncertain Job Skill Requirements}
\author[1]{Yulia Anoshkina\footnote{yulia.anoshkina@bwl.uni-kiel.de}}
\author[2]{Marc Goerigk\footnote{marc.goerigk@uni-siegen.de}}
\author[1]{Frank Meisel\footnote{meisel@bwl.uni-kiel.de}}

\affil[1]{School of Economics and Business, Christian-Albrechts-University of Kiel, Germany}
\affil[2]{Network and Data Science Management, University of Siegen, Germany}

\begin{document}
\maketitle	
\begin{abstract}
We consider a combined problem of teaming and scheduling of multi-skilled employees that have to perform jobs with uncertain qualification requirements. We propose two modeling approaches that generate solutions that are robust to possible data variations. Both approaches use variants of budgeted uncertainty, where deviations in qualification requirements are bounded by a constraint.
In the first approach, we aggregate uncertain constraints to ensure that the total number of job qualifications present at a job is not less than a worst-case value. We show that these values can be computed beforehand, resulting in a robust model with little additional complexity compared with the nominal model.
In our second approach, we bound the overall qualification deviation over all jobs. While this approach is more complex, we show that it is still possible to derive a compact problem formulation by using a linear programming formulation for the adversarial problem based on a dynamic program. 
The performance of both approaches is analyzed on a test bed of instances which were originally provided for a deterministic problem version. Our experiments show the effectiveness of the proposed approaches in the presence of data uncertainty and reveal the price and gain of robustness.
	\end{abstract}

	

\section{Introduction} \label{sec:Introduction}

This paper addresses a combined problem of routing and scheduling of multi-skilled workforce as it is faced by many service-oriented companies that provide installation, construction, maintenance or delivery services at customer locations. Each service job to conduct requires employees with different skill domains and at different levels of expertise. Therefore, teams of technicians have to be formed according to job qualification requirements that express the number of employees with specific skills and required experience in the corresponding domains.
In order to increase productivity and to decrease labor costs, companies may prefer to hire multi-skilled employees that  can be easily assigned to various jobs as required. This provides more flexibility and allows a company to focus on customer satisfaction. As teams may be capable of serving multiple jobs, optimal routing plans have to be found for the formed teams. From this, the investigated Routing and Scheduling Problem of Multi-Skilled Teams (\problem)
can be considered as an extension of the Vehicle Routing Problem (VRP). Due to its practical relevance, the {\problem}
has gained an increasing attention during the last decade and has been investigated extensively from different perspectives. In this paper, we demonstrate how the {\problem} can be solved in the presence of data uncertainty. 
In general, different sources for data variations are existent but we focus here on uncertainty of the qualification requirements of a job. 
Such qualification requirements are  usually derived from communication with customers. As a customer is not necessarily an expert in the corresponding field, the required skill types and levels of competence for executing a job  may be wrongly assessed by the customer when issuing the job. 
Also the company may misproject these requirements before having executed the job due to a lack of information. 
Moreover, the routing decisions can be affected by variations of travel times due to traffic conditions or by delays in job processing due to differing employee working speeds.
In light of these findings, it becomes important to have a robust planning approach that ensures solution reliability also in the presence of possible data variations. However, despite the substantial progress in the field of robust optimization, we are not aware of any approaches that have been so far presented for the formation of worker teams and their job-routing as is addressed in this paper.  

Recently, \cite{anoshkina2019technician} analyzed the deterministic version of the {\problem}
and the potential of decomposition techniques for reducing the complexity of the planning. Furthermore, {\problem}
was considered from a multi-period perspective by additionally emphasizing team consistency in \cite{anoshkina2020interday}. In the present study, we take the first step to deal with data variations in the context of scheduling of multi-skilled teams. Namely, we concentrate on developing a linear optimization framework incorporating \textit{demand uncertainty} which we define as a variation of job qualification requirements. 
Our contribution is then threefold: (i) We propose a first robust model formulation based on the concept of {budgeted} uncertainty sets as  proposed by \cite{bertsimas2003robust}, where the uncertainty affects each job independently. (ii) We propose a second robust model formulation, where the uncertainty is restricted by a global constraint over all jobs. (iii) We test extensively the model performance under the two different robustness strategies by analyzing the impact of the uncertain demand on the  feasibility and quality of solution.

The outline of this paper is as follows. In Section~\ref{sec:Literature}, we review the relevant literature on workforce teaming and scheduling. We also discuss important robustness concepts that constitute the foundation of our later investigation. In Section~\ref{sec: model}, we present a mathematical formulation of the deterministic problem version. In Section~\ref{sec: robust formulation}, we develop two robust optimization models based on different budgeted uncertainty sets. Section~\ref{sec:Computational study} presents experimental results and analyzes the performance of the two robustness strategies. Finally, Section~\ref{sec:Conclusions} concludes the paper and outlines future research.

\section{Literature} \label{sec:Literature}

The combined problem of teaming and scheduling of a multi-skilled workforce 
was first addressed in the works of \cite{estellon2009high}, \cite{hurkens2009incorporating},  \cite{cordeau2010scheduling}, \cite{hashimoto2011grasp}, and \cite{firat2012improved}. The initial focus of the research has been on scheduling aspects. Specifically, it was considered how multi-skilled employees can be grouped into teams and assigned to a set of jobs where jobs require multiple skills at different competence levels. 
Following these initial contributions, an increasing number of extensions has been presented. 
For instance, such features as
routing of teams \citep{kovacs2012adaptive}, multi-period planning \citep{zamorano2017branch}, employee preferences for performing a specific job \citep{firat2016branch}, alternative heuristic solution methods \citep{khalfay2017applying}, decomposition techniques \citep{anoshkina2019technician} as well as team consistency and rescheduling \citep{anoshkina2020interday} have been studied.
A more detailed description of these studies is provided by \cite{anoshkina2020interday}. 
Throughout, the authors assumed a deterministic setting where all input data is completely known with certainty.

However, real-world situations typically involve data uncertainty. 
Therefore, considerable research has been conducted in developing robust programs that 
find solutions which perform well despite variations in the input data. 
A substantial progress in the theory of robust optimization has been achieved with concepts presented by \cite{ben1999robust} and \cite{bertsimas2004price}. More precisely, \cite{ben1999robust} showed that robust counterparts of linear programs with ellipsoidal uncertainty set are computationally tractable and can be solved as conic quadratic problems.
\cite{bertsimas2004price} developed the concept of budgeted uncertainty that enabled to reformulate non-linear robust constraints as linear functions. 
For general surveys on robust optimization, we refer the interested reader to \cite{gabrel2014recent,goerigk2016algorithm,buchheim2018robust}.

The mentioned approaches have opened an avenue for research in many optimization areas. For instance, \cite{sungur2008robust} addressed a VRP with uncertain customer demand and proposed three robust formulations based on convex, box and ellipsoidal uncertainty sets. Following the idea of \cite{bertsimas2003robust}, \cite{ordonez2010robust} presented a robust formulation of the VRP incorporating two additional sources of uncertainty occurring in travel time and travel costs. 
This line of research was continued for a number of problem extensions. 
For instance, \cite{lee2012robust} addressed a VRP with deadlines and uncertainty arising in travel times and customer demand. 
\cite{han2014robust} combined stochastic programming with robust optimization for the solution of a VRP with uncertain travel times where penalties are imposed if travel time exceeds a preset time limit.
Demand uncertainty was also studied by \cite{cao2014open} in the context of open VRPs where vehicles do not necessarily return to the depot after delivering goods.
\cite{chen2016robust} analyzed a routing problem arising in road maintenance, in which each part of a road network has to be monitored by a service vehicle. Thereby, service times are subject to uncertainty due to various factors like road conditions or accidents. 
\cite{de2019robust} investigated a VRP with time windows (VRPTW) and multiple delivery men where a specific number of workers is required to execute deliveries and customer demand becomes known only when a vehicle arrives at a customer location. 
A further VRPTW model with both demand and travel time uncertainty was provided by \cite{munari2019robust} where the authors used a  two-index vehicle flow formulation. The main advantage of this formulation is that the robust counterpart can be derived directly from the underlying deterministic model and, thus, does not require additional constraints associated with uncertain parameters.

Compared with the large number of studies addressing robust VRPs, relatively few papers have been published on robust personnel scheduling. \cite{carello2014cardinality} developed a robust optimization model for a home health care problem with demand uncertainty and continuity of care. Continuity of care means that all services required by a patient are provided by the same specialist over a long period. The demand is considered uncertain due to possible variations in the physical conditions of patients. \cite{nguyen2016robust} proposed a nonlinear 
mixed-integer programming
formulation for taking into account uncertainty in nurse availability. 
\cite{souyris2013robust} examined the problem of dispatching technicians under stochastic service times. Specifically, the authors developed two different solution concepts distinguishing between processing time uncertainty related to customers or to  technicians.
Finally, we are aware of only one robust optimization approach dealing with scheduling of multi-skilled employees where the workforce demand is subject to uncertainty. At the example of a service industry company, \cite{henao2016multiskilling} investigate how multi-skilled employees can be effectively distributed between departments over a planning horizon of one week. Thereby, the problem also incorporates decisions about training of employees specialized only in one domain. The goal is to minimize staff training, shortage and surplus costs. Operations management for the assignment of jobs to teams and routing decisions, as is considered in our paper, are out of scope of their study.

As results generated by robust programs can deviate significantly from deterministic solutions, a further stream of research focuses on methods and algorithms that reduce this so-called price of robustness. Complementing the work of \cite{bertsimas2003robust}, \cite{poss2013robust} presented the concept of variable budgeted uncertainty, where dualization techniques are applied to more general uncertainty polytopes. It was shown by experiments that the proposed approach can yield better results and reduces the price of robustness by 18\%. 
Furthermore, the robust optimization methodology was extended by a class of two-stage robust optimization concepts, see e.g.  \cite{ben2004adjustable,hanasusanto2015k,liebchen2009concept,adjiashvili2015bulk,buchheim2017min}. These approaches consider problems where decisions can be taken sequentially. Therefore, a subset of decisions is implemented before the specific data realization becomes known whereas the remaining decisions can be taken after the uncertainty has resolved. For a general survey on these approaches, we refer to \cite{yanikouglu2019survey}.

To the best of our knowledge none of the mentioned concepts has been applied so far to the {\problem}. To close this gap, we develop two alternative robust optimization models that incorporate demand uncertainty in the context of routing and scheduling of multi-skilled teams.

\section{Deterministic Model (\textbf{DM})} \label{sec: model}\label{sec: deterministic model }

The deterministic (nominal) version of the multi-skilled workforce routing and scheduling problem can be described as follows. We are given a set of employees $M$ and a set of jobs $J$ as well as an extension of this set as $J^{0}=\{0\}\cup J=\{0,1,\ldots,|J|\}$ where 0 refers to a depot. Each job $j \in J$ is characterized by a service requirement $r_{jkl}$ that gives the number of employees  with qualification in skill $k \in K$ and experience level $l \in L$ required for performing job $j$. Here, $K$ denotes the set of skill domains and $L$ the set of experience levels. 
The competences of employee $m \in M$ are described by a binary matrix $q_{mkl}$ where an element takes value 1 if the employee is qualified in skill $k \in K$ at level $l \in L$ and 0 otherwise. As each job can require more than one employee, 
employees have to be grouped into teams in order to meet a job's qualification requirements $r_{jkl}$ for all $k \in K$ and $l \in L$. 
The maximal number of teams $T$ to build is specified by the minimum of the number of employees and the number of jobs considered in a problem instance. More precisely, if we consider a problem with 10 employees and 5 jobs, at most $T=min\{10,5\}$ teams are required (or can be built). 
Note that each job $j \in J$ has to be carried out by exactly one team, whereas a team might perform several jobs one after the other. Thereby, the completion time of each job cannot be later than a maximal working time $e_{max}$ that is given for each team.
Further, all services associated with job $j$ are provided at the customer's location. To each pair of jobs $(i,j) \in J^{0}\times J^{0}$, we thus assign  a travel time $d_{ij}$ that is needed by a team to go from $i$ to $j$. Additionally, each job $j$  has a processing time $p_{j}$ that indicates the amount of time that a team has to stay at customer location $j$. Here, we assume that $p_{j}$ is given and constant, i.e. $p_{j}$ does not depend on the team composition or working environment. 

A corresponding deterministic model has been provided by \cite{anoshkina2020interday}. Here we present a slightly modified formulation that constitutes the foundation of our robust approach. 
The formulation uses the following decision variables.
The binary decision variable $x_{mt}$ indicates if employee $m$ is assigned to team $t$ or not.
The routes of teams are denoted by binary decision variables $z_{tij}$ that define if team $t$ travels directly from job $i$ to job $j$  or not.  The continuous scheduling variable $s_{tj}$ specifies the start time of job $j$ by team $t$. Similar, $f_{tj}$ denotes the completion time of job $j$ executed by team $t$. 
Using the introduced notation the deterministic model is formulated as follows.
\begin{linenomath}\begin{align}
\text{maximize:} \ & \alpha\cdot \sum_{t \in T}\sum_{j \in J^{0}}\sum_{j \in J} z_{tij} - \beta\cdot \sum_{t \in T}\sum_{j \in J} f_{tj} \label{obj}\\
\text{subject to:}\notag\\ 
\sum_{t \in T} x_{mt} & \leq 1 && \forall  m \in M\label{x<1}\\
\sum_{m \in M}x_{mt} \cdot q_{mkl} & \geq r_{jkl} \cdot \sum_{i \in J^{0}}z_{tij} && \forall  j \in J, k \in K, l \in L, t \in T\label{quali}\\
\sum_{j\in J} z_{t0j} &\leq 1 && \forall t \in T\label{depot=1}\\
\sum_{t \in T} \sum_{i \in J^{0}} z_{tij}& \leq 1 && \forall j \in J\label{(i,j)<1}\\
\sum_{i \in J^{0}} z_{tij} &= \sum_{i \in J^{0}} z_{tji} && \forall  j \in J^{0}, t \in T\label{no segmentation}\\
f_{ti}  + d_{ij}    &\leq s_{tj} + \mathcal{M} \cdot\left(1-z_{tij}\right) && \forall i \in J^{0},  j \in J,  t \in T \label{stj}\\
s_{tj} + p_{j}   & \leq f_{tj} + \mathcal{M} \cdot \left(1 - \sum_{i \in J^{0}}z_{tij}\right) && \forall j \in J, t \in T \label{ftj}\\
f_{tj}&\leq e_{max}  && \forall j \in J, t \in T \label{emax}\\
s_{tj}, f_{tj} & \geq 0 && \forall j \in J^{0}, t \in T\label{s,f, domain}\\
x_{mt},  z_{tij} & \in \lbrace 0,1\rbrace &&\forall i,j \in J^{0}, m \in M, t \in T\label{z, x domain}
\end{align}\end{linenomath}
The main goal of the model is to maximize the service level, which we define as the number of performed jobs. 
Minimization of the total job completion time is considered as a subordinate objective. 
Weights $\alpha$ and $\beta$  are used for expressing different priorities of these two objectives.  Constraints~\eqref{x<1} forbid to assign an employee to more than one team. 
Constraints~\eqref{quali} guarantee that each formed team $t$ has an appropriate qualification for processing all jobs that are assigned to this team. 
Constraints~\eqref{depot=1} indicate that each active team starts from the depot. 
Constraints~\eqref{(i,j)<1} impose that each job can be served by at most one team. 
Constraints~\eqref{no segmentation} are the flow balancing constraints.
Constraints~\eqref{stj}-\eqref{ftj} define the start and completion times of job $j$ performed by team $t$.
Here, 
$\mathcal{M}$ denotes a sufficiently large positive value. 
Note that Constraints~\eqref{stj}-\eqref{ftj} also prevent subtours in the solution. 
Constraints~\eqref{emax} bound the longest working time for all teams.
Constraints~\eqref{s,f, domain}-\eqref{z, x domain} specify the domains of decision variables.

\section{Robust Formulations with Uncertain Job Qualification Requirements}\label{sec: robust formulation}

We present two robust problem formulations, which are based on different models to treat uncertainty in skill demand. For the ease of notation, we denote the two models by \textbf{RM1} (first robust model) and \textbf{RM2} (second robust model), respectively. 

\subsection{First Robust Model (\textbf{RM1})}

Our aim is to generate solutions that are insensitive to demand deviations. By demand deviation, we understand the variation of skill vectors in a job requirement matrix  $r_{jkl}$. 
In the deterministic model presented in Section~\ref{sec: model}, only Constraints~\eqref{quali} are affected by the variation of job skill requirements $r_{jkl}$. Note that only one element of $r_{jkl}$ is examined in each qualification constraint. 
In the standard technique for robust optimization \citep{bertsimas2004price}, we require constraint-wise uncertainty instead. Otherwise, we would hedge against the worst case in each parameter, which results in overly conservative solutions.
To avoid this conservatism, we follow the approach of \cite{bohle2010robust} and \cite{henao2016multiskilling} and extend the original deterministic model by redundant constraints expressing the aggregated qualification requirement, which we compute as the sum of technicians in all skill domains required on all levels of competence:  
\begin{linenomath}\begin{equation}
\sum_{m \in M}\sum_{k \in K}\sum_{l \in L}x_{mt} \cdot q_{mkl}  \geq \sum_{k \in K}\sum_{l \in L}r_{jkl}\cdot \sum_{i \in J^{0}}z_{tij}  \hspace{20 mm} \forall  j \in J, t \in T\label{quali reform}\\
\end{equation}\end{linenomath}
We model the uncertain demand $\tilde{r}_{jkl}$ for all $j \in J$ as an independent, random variable bounded on the interval 
$\tilde{r}_{jkl} \in [ r_{jkl}, r_{jkl} + \hat{r}_{jkl}]$, where $r_{jkl}$ denotes the nominal value and  $\hat{r}_{jkl}$ the maximal deviation allowed for $r_{jkl}$. 
For each random variable $\tilde{r}_{jkl}$, we define a level of variability $\zeta_{jkl}$ ranging within $[0,1]$. 
From this, the skill requirement variation is formulated as follows: 
\begin{linenomath}\begin{equation}
\tilde{r}_{jkl} = r_{jkl} +  \hat{r}_{jkl} \cdot \zeta_{jkl} \hspace{59 mm} \forall j \in J, k \in K, l \in L
\end{equation}\end{linenomath}
Furthermore, we assume that any skill and any qualification level can be exposed to uncertainty. The level of uncertainty for a job $j$ is controlled by parameter $\Gamma_{j}\in\mathbb{N}$ that presets the maximum skill and experience deviation allowed for this job. More precisely, 
$\Gamma_{j}$ represents an upper bound on the sum of skill and experience deviation weights $\zeta_{jkl}$ over all skill domains $k \in K$ and levels $l \in L$.  
From this, $\Gamma_{j}$ serves to adjust the robustness of the solution against the level of conservatism of a decision maker \citep{bertsimas2003robust}. For instance, if $\Gamma_{j}=0$, a decision maker assumes that no element of $r_{jkl}$ is likely to change. This corresponds to a risk seeking attitude where no protection against demand uncertainty is incorporated in the planning. 
In contrast, $\Gamma_{j}=\vert K\vert \cdot \vert L \vert$ assumes that all elements of $r_{jkl}$ are subject to uncertainty,  which corresponds to a very risk averse decision maker. This guarantees the maximal level of protection against all possible variations but, at the same time, results in the most conservative solution. Based on the previous notation, the uncertainty set for each job $U_{j}^{\Gamma}$  is defined as follows:
\begin{linenomath}\begin{equation*}
U_{j}^{\Gamma}=\left\{\tilde{r}_j \in \mathbb{R}^{\vert K\vert \cdot \vert L\vert}\, |\,  \tilde{r}_{jkl} = r_{jkl} +\hat{r}_{jkl} \cdot \zeta_{jkl}, \ 0 \leq \zeta_{jkl} \leq 1 \ \forall  k\in K, l\in L, \ \sum_{k \in K}\sum_{l \in L} \zeta_{jkl} \leq \Gamma_{j} \right\}
\end{equation*} \end{linenomath}
The aim of the robust model is to find solutions that remain feasible for all possible qualification requirements $\tilde{r}_j\in U^\Gamma_j$ for each job $j$.

Using the uncertainty set $U_{j}^{\Gamma}$, the robust counterpart of \eqref{quali reform} can be formulated as
\begin{linenomath}\begin{equation}
\sum_{m \in M}\sum_{k \in K}\sum_{l \in L}x_{mt} \cdot q_{mkl}  \geq \sum_{k \in K}\sum_{l \in L} 
\tilde{r}_j
\cdot \sum_{i \in J^{0}}z_{tij}  \hspace{10 mm} \hspace{9 mm}\forall  j \in J,\tilde{r}_j \in U_{j}^{\Gamma}, t \in T\label{quali robust reform}
\end{equation}\end{linenomath}
Formulation~\eqref{quali robust reform} is intractable in its current form since it contains an infinite number of constraints for all realizations of the continuous parameters $\zeta_{jkl}$  within the uncertainty set $U_{j}^{\Gamma}$.
To approach this issue, note that $\sum_{i \in J^{0}}z_{tij}$ is either 0 or 1. Hence, there are only two cases we need to consider: If $\sum_{i \in J^{0}}z_{tij}=0$, constraint~\eqref{quali reform} is always fulfilled. If $\sum_{i \in J^{0}}z_{tij}=1$, then we need to calculate the maximum value that $\sum_{k \in K}\sum_{l \in L}\tilde{r}_{jkl}$ can possibly take. Denoting this value by $\bar{r}_j$, constraint~\eqref{quali robust reform} becomes
\begin{linenomath}\begin{equation}\label{robcon}
\sum_{m \in M}\sum_{k \in K}\sum_{l \in L}x_{mt} \cdot q_{mkl}  \geq \bar{r}_j \cdot \sum_{i \in J^{0}}z_{tij}  \hspace{20 mm} \forall  j \in J, t \in T
\end{equation}\end{linenomath}
To calculate $\bar{r}_j$, we need to solve the problem
\begin{linenomath}\begin{align*}
\bar{r}_j = &\max \left\{ \sum_{k\in K} \sum_{l\in L} (r_{jkl} + \hat{r}_{jkl}\cdot \zeta_{jkl}) : \sum_{k\in K} \sum_{l\in L} \zeta_{jkl} \le \Gamma_j,\ 0 \le \zeta_{jkl} \le 1\  \forall k\in K, l\in L \right\} \\
=&\sum_{k\in K} \sum_{l\in L} r_{jkl} + \max \left\{ \hat{r}_{jkl}\cdot \zeta_{jkl} : \sum_{k\in K} \sum_{l\in L} \zeta_{jkl} \le \Gamma_j,\ 0 \le \zeta_{jkl} \le 1 \ \forall k\in K, l\in L \right\}
\end{align*}\end{linenomath}
Calculating this value can be done by sorting the vector $\hat{r}_{jkl}$ in descending order, and then summing up the $\Gamma_j$ many largest values. Then, the robust counterpart of the nominal model becomes
\begin{linenomath}\begin{flalign}
\text{maximize: } & \alpha\cdot \sum_{t \in T}\sum_{i \in J^{0}}\sum_{j \in J} z_{tij} - \beta\cdot \sum_{t \in T}\sum_{j \in J} f_{tj}\\
\text{subject to: } & \eqref{x<1} -\eqref{z, x domain} \ \text{and} \notag\\ 
& \sum_{m \in M}\sum_{k \in K}\sum_{l \in L}x_{mt} \cdot q_{mkl}  \geq \bar{r}_j \cdot \sum_{i \in J^{0}}z_{tij}  \hspace{20 mm} \forall  j \in J, t \in T \label{cons 17}
\end{flalign}\end{linenomath}
We further extend this model by measuring how much the required right-hand side constraint~\eqref{cons 17} is exceeded. This excess creates an additional benefit for the objective function, i.e., we reward additional robustness in the solution with some weight $\mu$. To this end, we introduce a new variable $\rho_{jt}$ that measures the slack of the right-hand side. 
This yields model \textbf{RM1}:
\begin{linenomath}\begin{flalign}
\text{maximize: } & \alpha\cdot \sum_{t \in T}\sum_{i \in J^{0}}\sum_{j \in J}  z_{tij} - \beta\cdot \sum_{t \in T}\sum_{j \in J} f_{tj} + \mu \cdot \sum_{t \in T} \sum_{j\in J} \rho_{jt} \\
\text{subject to: } & \eqref{x<1} -\eqref{z, x domain} \ \text{and} \notag\\ 
&  \sum_{m \in M}\sum_{k \in K}\sum_{l \in L}x_{mt} \cdot q_{mkl} \ge \bar{r}_j \cdot\sum_{i \in J^{0}}z_{tij} + \rho_{jt}  &  \forall  j \in J, t \in T  \\
& \rho_{jt} \le \mathcal{M}\cdot \sum_{i \in J^{0}}z_{tij} & \forall  j \in J, t \in T \label{adjconst} \\
& \rho_{jt} \ge 0 & \forall j\in J, t \in T
\end{flalign}\end{linenomath}
The additional constraint~\eqref{adjconst} is required to ensure that the excess is only taken into account if the job $j$ is actually performed by team $t$.

\subsection{Second Robust Model (\textbf{RM2})}

In the previous model, uncertainty sets were applied job-wise, which makes it possible to find a robust counterpart with little computational overhead. It has the drawback that solutions may still become overly conservative, as worst-case scenarios are assumed for each job separately. Furthermore, the aggregation of constraints means that it is ignored with what skills we hedge against uncertainty, as long as the total number of skills present is sufficient. We now follow a more nuanced approach to model uncertainty, which avoids both problems.

Consider the skill requirement $r_{jkl}$ for job $j$, skill $k$, level $l$. Let us assume we build a team that reaches a qualification level $\sum_{m \in M}x_{mt}\cdot q_{mkl}$. The buffer is then defined as $b_{jkl} = \sum_{m \in M}x_{mt}\cdot q_{mkl} - r_{jkl}$.

Let us assume there is an adversary who tries to find a scenario to disrupt as many jobs as possible. The adversary can increase the required skill level $r_{jkl}$ under the following conditions: increasing $r_{jkl}$ by one unit has some cost $c_{jkl}$, which reflects that higher level skills are less likely than lower level skills ($c_{jkl}$ increases with $l$) and that it should be more expensive to increase the demand of skills $k$ that are less likely to be relevant as judged by expert knowledge. The adversary has a global budget $\Gamma$ he can use for skill requirement increases.
A job is disrupted if the requirements in one skill and level are not met.

Given a fixed team and schedule, we hence want to solve the following adversary problem:
\begin{linenomath}
\begin{align}
\text{maximize: }\ &\sum_{j\in J} \left(\sum_{t\in T} \sum_{i\in J^0} z_{tij}\right)\cdot \zeta_j \\
\text{subject to: } & \zeta_j \le \sum_{k\in K} \sum_{l\in L} \zeta_{jkl} & \forall j\in J \label{advcon1} \\
& \sum_{j\in J} \sum_{k \in K} \sum_{l \in L} (b_{jkl}+1)\cdot c_{jkl}\cdot  \zeta_{jkl} \le \Gamma \label{advcon2}\\
& \zeta_j \in \{0,1\} & \forall j\in J\\
& \zeta_{jkl} \in \{0,1\} & \forall j\in J, k\in K, l\in L
\end{align}
\end{linenomath}
Here, binary variable $\zeta_{jkl}$ indicates if job $j$ is prevented by increasing the requirements in skill $k$ at level $l$. Binary variable $\zeta_j$ indicates if job $j$ is prevented overall. 
Using constriant~\eqref{advcon1}, we enforce that $\zeta_j$ can only be active if at least one of the $\zeta_{jkl}$ variables is active as well. Constraint~\eqref{advcon2} ensures that the total budget is restrcited to $\Gamma$, where the costs on the left-hand side correspond to the required investment to make a job infeasible.
Note that in an optimal solution, one would not increase multiple variables $\zeta_{jkl}$ for the same job $j$, but only choose the cheapest possibility. As the buffers $b_{jkl}$ depend on the assignment $x_{mt}$, we do not remove these variables from the problem.

Unfortunately, it is not possible to relax this formulation of the adversarial problem without changing its objective value. This means that a compact robust formulation cannot be obtained by simply dualizing the linear relaxation of the adversarial problem. In the following, we show that a compact formulation can still be obtained by using a dynamic programming formulation.

Let us denote by $F(j,\gamma)$ the maximum number of jobs from $j'\in\{0,\ldots,j\}$ that can be interrupted with a budget $\gamma\in\{0,\ldots,\Gamma\}=:\Gamma^0$. We have $F(0,\gamma)=0$ for all $\gamma\in\Gamma^0$, and the recursion
\begin{linenomath}\[ F(j,\gamma) = \max\left\{ F(j-1,\gamma), 1+F(j-1,\gamma-\min_{k,l} (b_{jkl}+1)\cdot c_{jkl}) \right\} \]\end{linenomath}
The value $F(|J|,\Gamma)$ is then equal to the objective value of the adversarial problem. We can also see this dynamic program as a longest path problem. We define a set of nodes $V=J^0\times\Gamma^0$
and arcs $A = \{ (j,\gamma,j',\gamma')\in V\times V : j' > j, \gamma' > \gamma\}$.
The adversary problem is then equivalent to solving
\begin{linenomath}
\begin{align}
\text{maximize: }\ & \sum_{a=(j,\gamma,j',\gamma')\in A} \left(\sum_{t\in T} \sum_{i\in J^0} z_{tij'}\right)\cdot  \tilde{c}_a p_a \\
\text{subject to: } & p \text{ is a path from } (0,0) \text{ to } (|J|,\Gamma)
\end{align}
\end{linenomath}
where
\begin{linenomath}\[
\tilde{c}_{j\gamma,j'\gamma'} = \begin{cases}
1 & \text{ if } \exists k,l: (b_{j'kl}+1)\cdot c_{j'kl} \le \gamma'-\gamma \\
0 & \text{ else} \end{cases}
\]\end{linenomath}
Dualizing this gives the model
\begin{linenomath}
\begin{align}
\hspace{-1 cm} \text{minimize: } & u_{|J|,\Gamma} \\
\hspace{-1 cm} \text{subject to: } & u_{00} = 0 \\
\hspace{-1 cm} & u_{j'\gamma'} \ge u_{j\gamma} + \left(\sum_{t\in T} \sum_{i\in J^0} z_{tij'}\right)\cdot \tilde{c}_{j,\gamma,j',\gamma'} & \forall (j,\gamma,j',\gamma') \in V\times V: j'>j, \gamma'>\gamma 
\end{align}
\end{linenomath}
Combining this dual adversarial model with the deterministic formulation, we obtain the following compact formulation \textbf{RM2} for the robust RSPMST:
\begin{linenomath}
\begin{flalign}
\text{maximize:} \ & \alpha\cdot\sum_{t \in T}\sum_{j \in J^{0}} \sum_{j \in J}  z_{tij} 
- \nu\cdot u_{|J|\Gamma} + \mu \cdot \sum_{j \in J}\sum_{k \in K}\sum_{l \in L}\sum_{t \in T} \rho_{jklt}
- \beta\cdot \sum_{t \in T}\sum_{j \in J} f_{tj} \hspace{-10cm} \label{comp1}\\
\text{subject to: } & \eqref{x<1}, \eqref{depot=1}-\eqref{emax} \text{ and}\notag \\  
\sum_{m \in M}x_{mt} \cdot q_{mkl} & \geq  r_{jkl} \cdot \sum_{i \in J^{0}}z_{tij} + \rho_{jklt}  && \forall  j \in J, k \in K, l \in L, t \in T \label{comp3}\\
\rho_{jklt} &\le \mathcal{M}\cdot \sum_{i\in J^0} z_{tij} && \forall j \in J, k \in K, l \in L, t \in T \label{comp4}\\
u_{00} &= 0 \label{comp7} \\
u_{j'\gamma'} &\ge u_{j\gamma} + v_{j',\gamma'-\gamma} && \forall j'\in J,j\in J^0:j'>j, \gamma',\gamma\in\Gamma^0 : \gamma'\ge \gamma  \label{comp8}\\
\mathcal{M}\cdot(v_{j\gamma} + &(1 - \sum_{i \in J^{0}} z_{tij}))  \geq \notag\\
\gamma -(\sum_{m \in M} & q_{mkl}\cdot x_{mt} - r_{jkl} +1) \cdot c_{jkl} + 1 && \forall \gamma \in \Gamma^{0}, j \in J, k\in K, l\in L, t \in T\label{comp15}\\
s_{tj}, f_{tj}, \rho_{jklt} & \geq 0 && \forall j \in J^{0}, t \in T, k\in K, l\in L \label{comp16}\\
x_{mt},  z_{tij} & \in \lbrace 0,1\rbrace &&\forall i,j \in J^{0}, m \in M, t \in T \label{comp17}\\
v_{j\gamma} &\in\{0,1\} && \forall j\in J, \gamma\in \Gamma^0  \label{comp18}\\
u_{j\gamma} &\geq 0 && \forall j\in J^0, \gamma\in \Gamma^0 \label{comp19}
\end{flalign}
\end{linenomath}
The objective function consists of four components. The first is to maximize the number of jobs that are taken on. This is reduced by the number of jobs that can be interrupted by the adversary, weighed with a factor $\nu$. With a factor $\nu$ slightly smaller than one, we ensure that it is better to plan a job and then have it canceled, than not planning the job at all. The third component is to maximize buffer sizes, similar to model \textbf{RM1}.
The last component is the travel time. 
Constraints \eqref{comp3}-\eqref{comp4} are modified qualification requirements. 
Constraints~\eqref{comp7}-\eqref{comp15} are used to calculate $u_{|J|\Gamma}$, the number of interrupted jobs. To this end, the binary variable $v_{j\gamma}$ is forced to be 1 if $\sum_{t \in T}\sum_{i\in J^0}  z_{tij}=1$ (i.e., the job is being taken) and $\gamma$ is sufficiently large to disrupt job $j$.

\section{Computational Study}\label{sec:Computational study}

In this section, we describe the results of a computational study that aims at comparing the performance of the models described in Sections~\ref{sec: deterministic model } and \ref{sec: robust formulation} where we explore the effect of robust planning on the scheduling decisions. Next, we describe our experimental setup followed by a presentation of the obtained results.

\subsection{Experimental Setup}

Our experiments are based on the 12 instance sets of \cite{anoshkina2019technician}. Each instance set contains 10 instances and is distinguished according to the number of jobs and available employees. The first set contains small instances  with 4 jobs and 4 employees each, while the last set includes large instances with 20 jobs and 20 employees each.
All instances are available online at \url{www.scm.bwl.uni-kiel.de/de/forschung/research-data}.

In order to estimate the extent to which the skill variations can impact the solution quality, we use employee and job qualification matrices with $\vert K \vert=3$ skills and $\vert L\vert=3$ skill levels. From this, the maximum possible scaled skill deviation for each job in \textbf{RM1} is $\Gamma_{j}=3\cdot 3=9$.  In our  experiments, we limit $\Gamma_{j}$ for all jobs to value $\Gamma_{j}=4$, which corresponds to a medium level of risk aversion.
In contrast, the uncertainty budget $\Gamma$ for \textbf{RM2} has to be defined individually for all instances and instance sets. Therefore, preliminary experiments were conducted to estimate $\Gamma$ manually.
For \textbf{RM1}, we set the maximal skill deviation for each job as $\hat{r}_{jk1}=2$, $\hat{r}_{jk2}=1$ and $\hat{r}_{jk3}=1$. To evaluate the skill deviation for \textbf{RM2}, we define the cost matrix $c_{jkl}$ randomly as follows. The cost of increasing the skill requirement at level $l=1$ are set to 1 and 2 with an equal probability. Thereby, $c_{jk1}=1$ means that the corresponding $r_{jkl}$ element is more likely to be changed.  Similar, $c_{jk2} \in \{3, 4\}$ and $c_{jk3} \in \{5, 6\}$. 
The maximal working time $e_{max}$ is set to 540 minutes for all instances and models.
Putting emphasis on the service quality, we use the following parameters for evaluating the objective functions: $\alpha=1$, $\beta=0.0001$, $\mu=0.01$ and $\nu=0.99$.

All tests have been run on an Intel(R) Core (TM) i7-8700 3.20 GHz with 32 GB of RAM. We used CPLEX 12.10 for solving the mixed-integer programming models using a runtime limit of 3600 seconds per instance.

\subsection{Price of Robustness}

The first experiment is conducted to test the extent to which the proposed linear models can be solved to optimality and to examine the effect of the robust planning on scheduling decisions. In particular, we analyze the so-called price of robustness indicating the extent to which the optimal robust solution differs from the non-robust deterministic solution. As performance measure, we consider the difference in the achieved service levels, which we associate with the number and the complexity of performed jobs.

Table~\ref{table cplex} reports average results for all instance sets and each modeling approach obtained by CPLEX. The first column of the table shows the problem size. The next five columns display results for the deterministic optimization  model from Section~\ref{sec: model}, where the reported values are averages for the solutions of 10 instances in the corresponding instance set. The first column $Z$ shows the number of performed jobs. The second column $C$ indicates the average complexity of performed jobs. We define the job complexity as the average required skill in all skill domains and at all levels of competence for those jobs that are processed in a solution, i.e., $C=\sum_{j \in J^{sol}}\sum_{k \in K}\sum_{l \in L}r_{jkl}/\vert J^{sol}\vert$, where $J^{sol} \subset J$ denotes the jobs of the solution.
Further, columns $T$ and $E$ specify the number of active teams in the route plans and the number of employees assigned to these teams. The next column $F$ gives the total job completion time.  
The corresponding results for the robust optimization models are presented in the middle and at the right of the table.

Based on Table~\ref{table cplex}, the following differences in the performance of the models can be observed. As expected, we see that \textbf{DM} generates many solutions with a higher service level than (column $Z$) \textbf{RM1} and \textbf{RM2}. This is because \textbf{DM} considers only nominal qualification requirements without taking risks of data variation into account. Thereby, we see that the number of performed jobs increases for instances with a larger number of available employees (\textit{$|J| < |M|$}). 
In general,\textbf{ RM2} is more conservative as the service level achieved under \textbf{RM2} is slightly lower for nearly all instances, compared to\textbf{ RM1}. 
Another aspect is the complexity of the performed jobs. Looking at columns $C$, we see a clear tendency for \textbf{RM1} and \textbf{RM2} to avoid an assignment of more challenging jobs. 
However, no direct correlation can be derived.
Comparing both models, we observe lower (see e.g. instances $4\times 4$, $4\times 8$ or $6\times12$) as well as higher complexity values (see e.g. instances $6\times 6$, $8\times 6$, $10\times 7$).

\captionsetup[table]{singlelinecheck=false, justification=raggedright, skip=0.2cm}
\begin{table}[t]
\centering
\fontsize{10pt}{10pt}
\selectfont
\captionof{table}{Performance metrics for price of robustness}\label{table cplex}
\setlength{\tabcolsep}{1.5 mm}
\begin{tabular} {crrrrrcrrrrrcrrrrrcrrrrr}
\hline
Instance & \multicolumn{5}{c}{\textbf{DM}} & & \multicolumn{5}{c}{\textbf{RM1}} & & \multicolumn{5}{c}{\textbf{RM2}}\\
\cline{2-6}
\cline{8-12}
\cline{14-18}
\textit{$|J|$ $\times$ $|M|$}  & $Z$ & $C$ & $T$ & $E$ & $F$  & & $Z$ & $C$ & $T$ & $E$ & $F$  & & $Z$ & $C$ & $T$ & $E$ & $F$ \\
\hline

4 x 4 & 2.1 & 5.1 & 1.8 & 3.6 & 692  & & 1.7 & 3.9 & 1.3 & 4.0 & 591  & & 1.5 & 3.7 & 1.0 & 4.0 & 514 \\

4 x 8 & 3.2 & 5.0 & 3.1 & 7.8 & 1116 & &  2.9 & 4.6 & 2.5 & 8.0 & 1055  & &  2.1 & 4.2 & 1.6 & 8.0 & 765 \\

6 x 6 & 3.7 & 4.6 & 3.1 & 5.5 & 1329 & & 2.8 & 4.0 & 1.8 & 6.0 & 1047  & & 2.7 & 4.1 & 1.8 & 6.0 & 1010 \\ 

6 x 12 & 4.8 & 5.0 & 4.7 & 11.0 & 1706  & & 4.6 & 4.9 & 3.6 & 12.0 & 1787  & & 3.9 & 4.5 & 2.9 & 12.0 & 1506 \\

8 x 6 & 4.4 & 4.8 & 3.3 & 5.6 & 1553 & & 3.3 & 4.0 & 2.2 & 6.0 & 1215  & & 3.1 & 4.5 & 1.8 & 6.0 & 1136 \\

8 x 12 & 6.2 & 5.1 & 5.7 & 11.5 & 2206  & & 5.5 & 4.8 & 4.1 & 12.0 & 2062  & & 5.0 & 4.6 & 3.5 & 12.0 & 1892 \\ 

10 x 7 & 5.1 & 4.9 & 3.8 & 6.7 & 1750  & & 3.8 & 4.5 & 2.4 & 7.0 & 1366 & & 4.2 & 4.7 & 2.4 & 7.0 & 1574 \\

10 x 13 & 7.7 & 5.5 & 6.4 & 12.5 & 2802  & & 6.5 & 4.7 & 4.6 & 13.0 & 2443  & & 6.1 & 5.2 & 4.1 & 12.9 & 2327 \\

15 x 8 & 6.7 & 4.9 & 4.6 & 8.0 & 2376 & & 4.8 & 3.7 & 3.0 & 8.0 & 1738  & & 4.1 & 5.1 & 2.4 & 7.9 & 1500 \\

15 x 15 & 10.3 & 5.4 & 7.6 & 14.6 & 3897  & & 8.2 & 4.8 & 5.0 & 15.0 & 3115 & & 6.4 & 5.3 & 4.0 & 14.9 & 2452 \\ 

20 x 10 & 9.0 & 4.7 & 6.0 & 9.9 & 3157  & & 6.8 & 4.3 & 3.5 & 10.0 & 2456 & & 6.6 & 4.3 & 4.2 & 10.0 & 2423 \\ 

20 x 20 & 13.7 & 5.5 & 10.4 & 19.9 & 5001  & & 10.9 & 4.6 & 6.8 & 20.0 & 4037  & & 8.3 & 4.9 & 5.3 & 19.8 & 3057\\ 

\hline
\end{tabular}
\end{table}

Considering the number of teams and employees used in the solutions (columns $T$ and $E$), we observe that \textbf{RM1} assigns employees to a consistently lower number of teams than \textbf{DM}.
This indicates that larger teams are created in order to guarantee a greater schedule reliability in the presence of possible data variations. 
Moreover, we observe a further decrease of $T$ when comparing \textbf{RM1} and \textbf{RM2} for the majority of instances.
Similar to \textbf{RM1}, all employees are involved. An exception to this are large-sized instances that could not be solved to optimality, see instances $10\times 13$ - $15 \times 15$ and $20\times 20$. This also demonstrates a similar tendency to increase the level of protection by increasing the team size.

A further examination shows that, compared to \textbf{DM}, \textbf{RM1} and \textbf{RM2} result in a lower total job completion time due to a lower number of performed jobs.  

Table~\ref{table gap, cpu} provides statistics for the consumed runtime expressed in seconds (column $CPU$) and the optimality gap in percent (columns $GAP$) reported by CPLEX after the runtime limit of 3600 seconds per instance. 
The optimality gap is computed as $\textit{GAP}=\text{(Objective - LB)}/\text{LB}$ where \glqq Objective\grqq \ denotes the value of the objective function achieved by the model and \glqq LB\grqq \ gives the lower bound value reported by CPLEX. Column $Opt.$ gives the number of instances solved to optimality in each instance set. 

\captionsetup[table]{singlelinecheck=false, justification=raggedright, skip=0.2cm}
\begin{table}[t]
\centering
\fontsize{10pt}{10pt}
\selectfont
\captionof{table}{Performance metrics for computation times}\label{table gap, cpu}
\setlength{\tabcolsep}{1.5 mm}
\begin{tabular} {crrrcrrrcrrr}
\hline
Instance & \multicolumn{3}{c}{\textbf{DM}} & & \multicolumn{3}{c}{\textbf{RM1}} & & \multicolumn{3}{c}{\textbf{RM2}} \\ 
\cline{2-4}
\cline{6-8}
\cline{10-12}
\textit{$|J|$ $\times$ $|M|$} & $CPU$  & $GAP$ & $Opt.$ & & $CPU$ & $GAP$ & $Opt.$ & & $CPU$ & $GAP$ & $Opt.$ \\ 
\hline

4 x 4 &  0.02 & 0 & 10 & &  0.02 & 0 & 10 & & 0.16 & 0 & 10 \\

4 x 8 &  0.02 & 0 & 10 & & 0.04 & 0 & 10 & & 0.46 & 0  & 10\\

6 x 6 & 0.15 & 0  & 10 & & 0.70 & 0 & 10 & & 7.06 & 0 & 10\\

6 x 12 & 0.15 & 0 & 10 & & 1.27 & 0 & 10 & & 143.37 & 0 & 10\\

8 x 6 & 0.64 & 0 & 10 & & 0.79 & 0 & 10 & & 425.68 & 5  & 9\\

8 x 12 & 128.82 & 0 & 10 & & 346.74 & 0 & 10 & & 1669.68 & 22  & 6\\

10 x 7 & 142.72 & 0 & 10 & & 237.35 & 0 & 10 & & 1468.76 & 17 & 7\\

10 x 13 & 1621.64 & 2 & 7 & & 2720.58 & 7 & 4 & & 3522.07 & 55 & 1 \\

15 x 8 &  2865.84 & 33 & 3 & & 2883.81 & 29 & 3 & & 3600.00 & 81 & 0\\

15 x 15 & 3600.00 & 29 & 0 & & 3600.00 & 24 & 0 & & 3600.00 & 73 & 0\\

20 x 10 & 3600.00 & 54 & 0 & & 3600.00 & 42 & 0 & & 3600.00 & 76  & 0\\

20 x 20 & 3600.00 & 32 & 0 & & 3600.00 & 23 & 0 & & 3600.00 & 74  &0\\

\hline
\end{tabular}
\end{table}

The obtained results show that the computational time increases with an increase of the instance size. Looking at column $Opt.$, we see that already small instances containing less that 10 jobs could not be solved to optimality within the preset runtime limit. 
According to Table 2, there is no substantial difference in the complexity of \textbf{DM} and \textbf{RM1}. \textbf{RM1} delivers almost the same number of optimal solutions as \textbf{DM} does. Thereby, \textbf{RM1} demonstrates only slightly lower $GAPs$ values compared to \textbf{DM}, see instances $15\times 8$ - $20\times 20$. 
In contrast, \textbf{RM2} requires a considerably higher computational effort due to a much larger number of variables and constraints. Already for instances of size $8\times 6$, we observe a positive average $GAP$ and considerably larger CPU times.

\subsection{Benefit of Robustness}\label{sec: benefit of robustness}

The next two experiments assess the effect of data changes on the solution feasibility. In other words, we test how many planned jobs the teams can actually perform when uncertain skill requirements realize in the schedule execution. 
For this purpose, we generate for each optimization approach $1,000$ demand scenarios per instance set. Thus, the results are averages over $S=10\cdot1,000=10,000$ scenarios. 
We start by generating scenarios of type \textbf{RM1} which are modeled with $\Gamma_{j}=3$, i.e. 
3 elements  are varied in the original qualification requirement matrix of each job. 
The obtained results are reported in Table~\ref{table scenarios RM1}. The first column shows the problem size.
Columns $A$ in each block give the average of the absolute number of performed jobs 
while columns $R$ indicate  the average relative proportion of processed jobs in all scenarios in percent. For a scenario $s$, we compute $R_{s}$ as $R_{s}=A_{s}/Z$ where $Z$ refers to the number of originally performed jobs for the corresponding model.  Columns $B$ show the relative frequency in percent with which each robust model outperforms \textbf{DM}. 
Finally, we report in columns $W$ the relative frequency with which the service level attained by each robust model is lower than the nominal one.

\captionsetup[table]{singlelinecheck=false, justification=raggedright, skip=0.2cm}
\begin{table}[t]
\centering
\fontsize{10pt}{10pt}
\selectfont
\captionof{table}{Scenarios of type \textbf{RM1}. Best average service level per row is highlighted in bold.}\label{table scenarios RM1}
\setlength{\tabcolsep}{2.5 mm}
\begin{tabular} {crrcrrrrcrrrr}
\hline
Instance & \multicolumn{2}{c}{\textbf{DM}} & & \multicolumn{4}{c}{\textbf{RM1}}  & & \multicolumn{4}{c}{\textbf{RM2}} \\
\cline{2-3}
\cline{5-8}
\cline{10-13}
\textit{$|J|$ $\times$ $|M|$} & $A$ & $R$  & & $A$ & $R$ & $B$ & $W$ & & $A$ & $R$ & $B$ & $W$\\
\hline

4 x 4 & 0.50 & 21.33 &  & \textbf{1.11} & 63.96 & 50.63 & 2.64 & & 1.10 & 71.36 & 50.69 & 2.23 \\

4 x 8 & 1.02 & 33.35 &  & 1.46 & 50.12 & 44.84 & 12.90 & & \textbf{1.73} & 86.75 & 57.12 & 4.29 \\

6 x 6 & 0.65 & 18.23 &  & 1.76 & 67.08 & 73.50 & 5.74 & & \textbf{1.88} & 73.12 & 82.37 & 0.50 \\

6 x 12 & 1.66 & 34.93 &  & 2.59 & 55.13 & 62.60 & 12.75 & & \textbf{3.07} & 80.78 & 81.17 & 1.67 \\

8 x 6 & 0.78 & 18.30 &  & 1.96 & 60.80 & 75.05 & 4.82 & & \textbf{2.12} & 70.78 & 85.06 & 1.82 \\

8 x 12 & 1.37 & 22.77 &  & 3.02 & 54.11 & 78.49 & 5.44 & & \textbf{3.13} & 61.72 & 76.56 & 5.79 \\

10 x 7 & 0.96 & 20.12 &  & \textbf{2.17} & 56.90 & 73.29 & 2.89 & & 2.00 & 51.69 & 68.05 & 4.28 \\

10 x 13 & 1.39 & 18.23 &  & 2.98 & 45.96 & 77.00 & 9.15 & & \textbf{3.32} & 58.04 & 88.67 & 2.12 \\

15 x 8 & 1.21 & 18.55 &  & \textbf{2.63} & 53.74 & 73.50 & 5.84 & & 2.51 & 66.51 & 73.10 & 8.87 \\

15 x 15 & 1.64 & 16.19 &  & \textbf{4.23} & 51.55 & 91.00 & 1.65 & & 3.75 & 60.75 & 88.58 & 3.58 \\

20 x 10 & 1.55 & 17.47 &  & \textbf{3.30} & 49.08 & 81.87 & 4.15 & & 2.80 & 48.29 & 69.10 & 15.03 \\

20 x 20 & 1.70 & 11.78 &  & 5.03 & 46.06 & 96.20 & 0.92 & & \textbf{5.17} & 72.34 & 90.94 & 3.04 \\

\hline
\end{tabular}
\end{table}

The results for \textbf{DM} show that a considerable number of job assignments becomes infeasible. 
In fact, the relative service level $R$ drops below 20\% for the most instances. 
This means that although the deterministic model inserts a lot of jobs in a solution, it finally fails to process these jobs due to uncertain job requirements and insufficiently qualified teams. 
This low reliability can be substantially moderated by \textbf{RM1}, which is confirmed by significantly higher $R$ values ranging between 45\% and 67\%. However, note that these values already lie below 100\% for $\Gamma_{j}=3$. This is because $\Gamma_{j}=4$ guarantees the solution feasibility only for the aggregated skill level but not for every single element of matrix $r_{jkl}$. This means that a solution can become infeasible also for skill deviation that is below the defined uncertainty budget $\Gamma_{j}=4$.
Moreover, for all instances, we observe significantly higher absolute numbers of still feasible job assignments (see column $A$). Looking at column $B$, we see that \textbf{RM1} outperforms \textbf{DM} in 45\% to 96\% of all scenarios while $W$ ranges between 1\% and 13\% only.
Although  \textbf{RM1} achieved a lower service level than \textbf{DM} in the first experiment (see Table~\ref{table cplex}), 
 it now performs clearly stronger under the uncertain problem setting.
Also \textbf{RM2} is superior to \textbf{DM}. In fact, \textbf{RM2} delivers better results in 50\% to 90\% of all scenarios, see column $B$. 
Moreover, compared to \textbf{RM1}, 
we even observe a higher absolute and relative service level achieved under \textbf{RM2} for instances $4\times 8$ to $8\times 12$, $10\times 13$ and $20\times 20$. 
However, solutions for larger instances are less immunized against this type of uncertainty as they are solved under lower $\Gamma$ compared to small- and medium-sized instances. This is due to a higher problem complexity and, thus, computational effort required by \textbf{RM2} that does not allow to further increase uncertainty budget within the preset time limit.

\captionsetup[table]{singlelinecheck=false, justification=raggedright, skip=0.2cm}
\begin{table}[t]
\centering
\fontsize{10pt}{10pt}
\selectfont
\captionof{table}{Scenarios of type \textbf{RM2}. Best average service level per row is highlighted in bold.}\label{table scenarios RM2}
\setlength{\tabcolsep}{2.5 mm}
\begin{tabular} {crrcrrrrcrrrr}
\hline
Instance & \multicolumn{2}{c}{\textbf{DM}} & & \multicolumn{4}{c}{\textbf{RM1}}  & & \multicolumn{4}{c}{\textbf{RM2}} \\
\cline{2-3}
\cline{5-8}
\cline{10-13}
\textit{$|J|$ $\times$ $|M|$} & $A$ & $R$  & & $A$ & $R$ & $B$ & $W$ & & $A$ & $R$ & $B$ & $W$ \\
\hline

4 x 4 & 0.60 & 26.60 &  & 1.04 & 64.84 & 46.03 & 7.31 & & \textbf{1.14} & 77.39 & 51.56 & 4.98 \\

4 x 8 & 1.05 & 33.75 &  & 1.47 & 50.97 & 43.40 & 14.44 & & \textbf{1.69} & 84.75 & 54.42 & 6.71 \\

6 x 6 & 0.97 & 26.86 &  & 1.85 & 69.62 & 64.66 & 9.90 & & \textbf{1.98} & 75.19 & 74.24 & 2.76 \\

6 x 12 & 1.57 & 33.96 &  & 2.29 & 48.88 & 53.63 & 15.30 & & \textbf{3.03} & 79.83 & 79.32 & 2.40 \\

8 x 6 & 1.04 & 26.61 &  & 1.86 & 57.15 & 59.77 & 11.32 & & \textbf{2.21} & 73.79 & 75.44 & 3.01 \\

8 x 12 & 1.50 & 24.47 &  & 2.94 & 52.89 & 75.57 & 5.02 & & \textbf{3.23} & 63.91 & 81.24 & 6.31 \\

10 x 7 & 1.28 & 25.68 &  & 2.24 & 58.91 & 65.55 & 9.87 & & \textbf{2.27} & 57.85 & 71.19 & 8.05 \\

10 x 13 & 1.93 & 25.52 &  & 3.30 & 50.13 & 66.67 & 13.58 & & \textbf{3.46} & 60.54 & 76.22 & 7.21 \\

15 x 8 & 1.55 & 23.55 &  & \textbf{2.79} & 57.27 & 67.83 & 10.29 & & 2.76 & 71.43 & 71.17 & 9.18 \\

15 x 15 & 2.45 & 23.92 &  & \textbf{4.64} & 56.48 & 87.19 & 3.41 & & 4.09 & 66.31 & 78.06 & 7.31 \\

20 x 10 & 2.04 & 22.77 &  & \textbf{3.93} & 58.26 & 83.78 & 4.38 & & 3.16 & 53.38 & 67.58 & 13.94 \\

20 x 20 & 2.58 & 18.40 &  & \textbf{5.71} & 52.03 & 93.57 & 1.95 & & 5.30 & 72.91 & 86.73 & 5.89 \\

\hline
\end{tabular}
\end{table}

To achieve a fair comparison between the two robust approaches, we conduct a second experiment to evaluate the  quality of solutions under scenarios of type \textbf{RM2}. For this purpose, we create further $10,000$ scenarios ($1,000$ for each instance set) with uncertainty budget $\Gamma=10\cdot \vert J \vert$.
Following definition in Section~4.2, the scenarios are modeled such that the uncertainty budget is bounded over all jobs. 
To simulate different skill realizations, the sequence in which the jobs are considered is defined randomly for each scenario. For each considered job, one element is changed in matrix $r_{jkl}$. Thereby, skill domain and the competence level are selected randomly with equal probability of $1/\vert K\vert$ and $1/\vert L\vert$. The process is continued until the budget is reached.

The obtained results are summarized in Table~\ref{table scenarios RM2}. 
Here, we see that the general trends  
are similar to  those in the previous experiment. 
Compared to both robust approaches, we observe a low solution feasibility for \textbf{DM} with $R$ values that lie under 30\% for the majority of instances. With \textbf{RM1} this ratio is increased to further 48\% to 69\%. 
Moreover, we observe an increase in $R$ by further  6\% to 34\% when comparing \textbf{RM1} and \textbf{RM2}.
In absolute terms, \textbf{RM2} outperforms \textbf{RM1} in 67\%  of cases. This holds especially for instances that could be solved to optimality. Whereas, for the four last instances, \textbf{RM1} is again superior to \textbf{RM2}.  
Analyzing the performance of each single approach under different uncertainty settings, we observe higher \textit{B} and lower \textit{W} values in the previous experiment for both \textbf{RM1} and \textbf{RM2}. This is explained by a higher uncertainty imposed by scenarios of type \textbf{RM1} and, thus, a lower number of feasible job assignments in \textbf{DM}-solution (see column \textit{A} in Table~\ref{table scenarios RM1} and Table~\ref{table scenarios RM2}).

\subsection{Variation of Uncertainty Budget}

\begin{figure}[t]
\centering
\includegraphics[width=1\textwidth]{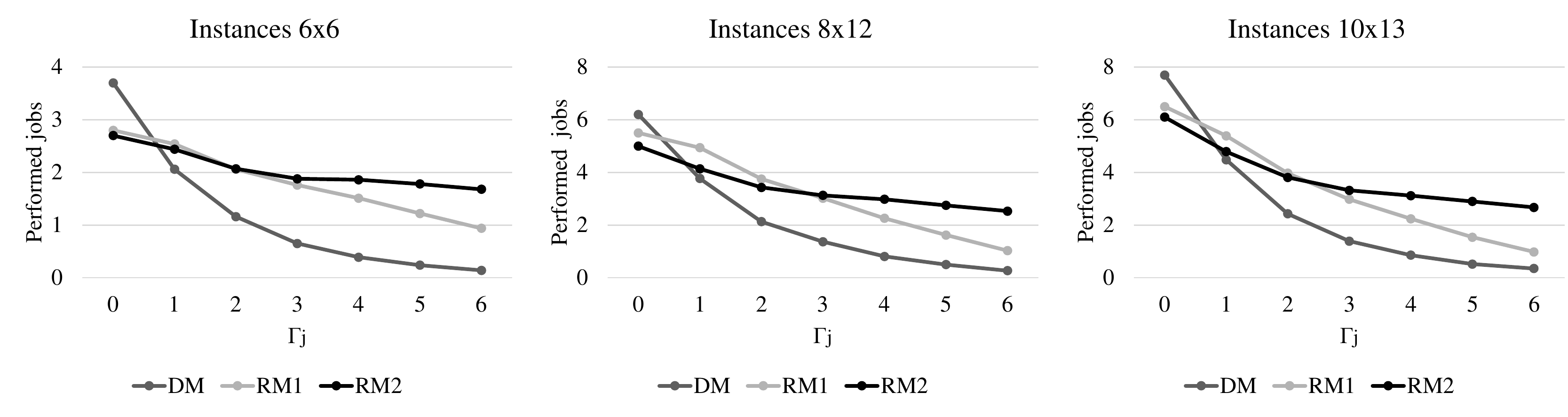}
\caption{Influence of uncertainty on the number of performed jobs with scenarios of type \textbf{RM1}}\label{Fig. gamma rm1}
\end{figure}

\begin{figure}[b]
\centering
\includegraphics[width=1\textwidth]{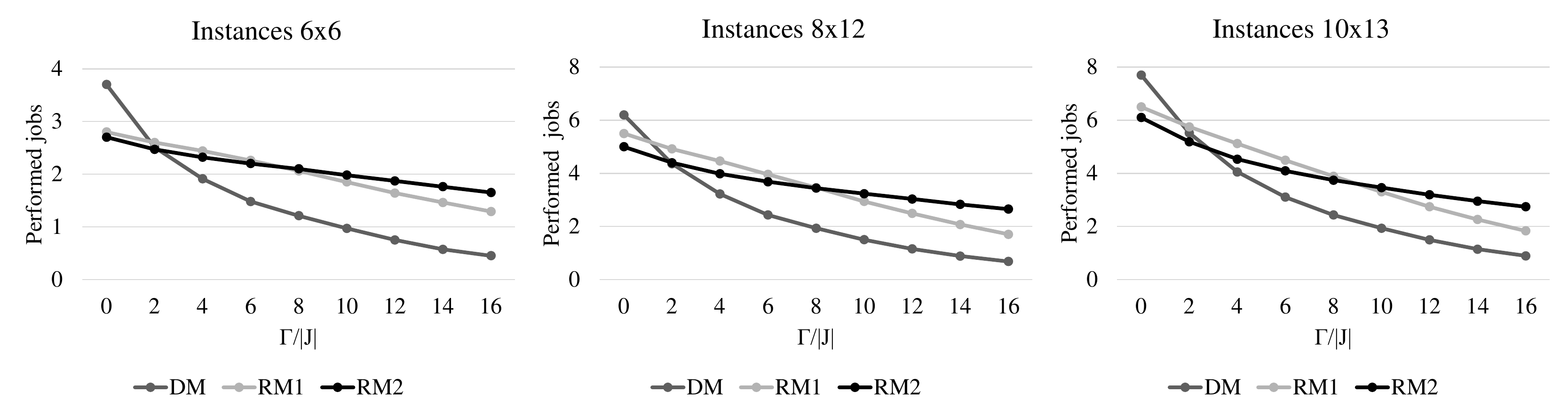}
\caption{Influence of uncertainty on the number of performed jobs with scenarios of type \textbf{RM2}}\label{Fig. gamma rm2}
\end{figure}

To give a more detailed understanding of the differences between the two robust planning approaches, we conduct a sensitivity analysis by varying the uncertainty budget for the generation of scenarios. Note that we do not recompute the model solutions in the simulation but evaluate feasibility of each scenario based on job assignments reported in Table~\ref{table cplex}. For scenarios of type \textbf{RM1}, $\Gamma_{j}$ is varied on the interval [0, 6]. Here, value 0 means that no skill deviations are considered, whereas value 6 means that six elements of the job requirement matrix can deviate from their nominal values simultaneously. For scenarios of type \textbf{RM2}, $\Gamma$ is varied on the interval $[0, 16\cdot \vert J\vert]$.
Figures~\ref{Fig. gamma rm1} and \ref{Fig. gamma rm2} show the impact of different parameter settings on the service level for three selected instance sets. Each plot relates to the average solutions of $10,000$ scenarios generated under different $\Gamma$-settings. Scenarios were generated according to the process described in the previous subsection. 

The results demonstrate that the service level is inversely correlated with the level of uncertainty. We can see a decline in the number of performed jobs with higher $\Gamma_{j}$ and  $\Gamma$ as the data becomes more and more uncertain. 
Thereby, marginal cost of robustness increase with an increase of the instance size. 
In general, solution quality in \textbf{RM1}-scenarios worsens at a strong rate. The service level drops to substantially lower values already in the middle of the examined interval. 
This is an expected outcome. As the uncertainty budget defined for \textbf{RM1} is aggregated job-wise, $\Gamma_{j}$-variations are more challenging and incur a higher price of robustness. 
Further, we observe that \textbf{RM1} achieves the same or even better service level at the interval $\Gamma_{j} \in [0, 3]$, whereas \textbf{RM2} generates a more robust solution for $\Gamma_{j} \in [4, 6]$, see Figure~\ref{Fig. gamma rm1}. To this end, the both proposed approaches significantly outperform \textbf{DM}.

A similar pattern emerges in Figure~\ref{Fig. gamma rm2}. Here, we also see that \textbf{RM1} is superior to \textbf{RM2} 
if the expected skill variation is relatively low. 
However, the break even point is now reached after only $\Gamma=8$. 
This is explained by the fact that \textbf{RM2} is more conservative with our choice of $\Gamma$. This confirm the results reported in Table~\ref{table cplex} where the service level achieved under \textbf{RM1} is consistently higher compared to \textbf{RM2}. From this, a higher level of data uncertainty is required by \textbf{RM2} to produce a significant impact on performance indicator.

To summarize our results, we can conclude that all proposed robust approaches can successfully handle the demand uncertainty. \textbf{RM2} provides a higher solution feasibility than \textbf{RM1} for our choice of $\Gamma$, but  tends to give the best performance on a wide range of levels of uncertainty, and even on scenarios that are generated in the style of \textbf{RM1}. The advantage of \textbf{RM1}, on the other hand, is its reduced computational complexity, which makes it possible to find robust solutions even on the largest instances.

\section{Conclusions and Future Research}\label{sec:Conclusions}

In this paper, we have investigated the problem of routing and scheduling multi-skilled teams under demand uncertainty where the variations of job qualification requirements are captured through uncertainty sets. For the solution of the problem, we have developed and analyzed two robust modeling approaches. Computational experiments showed that deviations in qualification requirements can have an extremely negative impact on the quality and the feasibility of the obtained solutions in a non-robust planning. This can be significantly moderated by the proposed robust approaches. The degree of solution robustness can be controlled not only by choosing uncertainty budget $\Gamma$ but also by choosing an appropriate method to model the uncertainty set. Specifically, we demonstrated that a higher protection against any demand variations is provided if the demand uncertainty can be distributed over the complete customer network where uncertainty cost are defined for each particular skill. 
Alternatively, uncertainty might be aggregated for each job separately. This allows to reach a reasonable compromise between the risk aversion and the achieved service level.

As this study represents a first step to incorporate uncertainty into scheduling of multi-skilled teams, there are still many promising avenues for future research. For instance, this study can be extended to other variants of uncertainty sets. 
In practice, the changes in job qualification requirements are often coupled with changes in job processing times. From this, it could be interesting to model interdependencies between these two parameters in the context of robust optimization. Finally, as the optimization models cannot be solved to optimality for large-scale problems, it would be worthwhile to develop further heuristic approaches.

\end{document}